\documentclass[a4paper,12pt]{article}
\makeatletter
\@addtoreset{footnote}{page}
\makeatother

\renewenvironment{enumerate}
{\begin{list}{\parbox{2em}{(\arabic{enumi})}}{
\usecounter{enumi}
\setlength{\topsep}{0ex}
\setlength{\itemindent}{1em}
\setlength{\leftmargin}{2em} %left indent
\setlength{\rightmargin}{1em}%right indent
\setlength{\labelsep}{1em}%sep between label and text
\setlength{\labelwidth}{3em}%
\setlength{\itemsep}{0.1em}%line break
\setlength{\parsep}{0em}%paragraph break
\setlength{\listparindent}{0em}%paragraph indent
}}{\end{list}}

\usepackage{amsthm}
\usepackage{amsmath,amssymb,latexsym,amsfonts,mathrsfs}

\renewcommand{\d}{{\rm d}} %differential

\newcommand{\dist}{\stackrel{{\rm d}}{=}}

\newcommand{\tend}[2]{\mathrel{\mathop{\longrightarrow}\limits^{#1}_{#2}}}

\newcommand{\supre}[2]{\mathrel{\mathop{\sup}\limits^{#1}_{#2}}}

\renewcommand{\hat}{\widehat}
\renewcommand{\tilde}{\widetilde}
\renewcommand{\bar}{\overline}

\newcommand{\absol}[1]{\left| #1 \right|} %absolute value
\newcommand{\rbra}[1]{\!\left( #1 \right)} %round brackets or parentheses
\newcommand{\cbra}[1]{\!\left\{ #1 \right\}} %curly brackets or braces
\newcommand{\sbra}[1]{\!\left[ #1 \right]} %square brackets or brackets
\newcommand{\bE}{\ensuremath{\mathbb{E}}}

\newcommand{\bN}{\ensuremath{\mathbb{N}}}

\newcommand{\bP}{\ensuremath{\mathbb{P}}}

\newcommand{\bR}{\ensuremath{\mathbb{R}}}

\newcommand{\cB}{\ensuremath{\mathcal{B}}}

\newcommand{\cF}{\ensuremath{\mathcal{F}}}

\newcommand{\cN}{\ensuremath{\mathcal{N}}}

\newcommand{\cP}{\ensuremath{\mathcal{P}}}

\newcommand{\cX}{\ensuremath{\mathcal{X}}}
\newcommand{\cY}{\ensuremath{\mathcal{Y}}}

\theoremstyle{plain}
\newtheorem{Thm}{Theorem}[section]

\newtheorem{Lem}[Thm]{Lemma}
\newtheorem{Prop}[Thm]{Proposition}
\newtheorem{Cor}[Thm]{Corollary}

\theoremstyle{definition}

\newtheorem{Def}[Thm]{Definition}

\newtheorem{Ex}[Thm]{Example}

\newcommand{\Proof}[2][Proof]{\begin{proof}[{#1}] #2 \end{proof}}

\numberwithin{equation}{section}

\makeatletter
\renewcommand\section{\@startsection {section}{1}{\z@}%
                                   {-3.5ex \@plus -1ex \@minus -.2ex}%
                                   {2.3ex \@plus.2ex}%
                                   {\normalfont\large\bf}}
\makeatother

\makeatletter
\renewcommand\subsection{\@startsection {subsection}{1}{\z@}%
                                   {-3.5ex \@plus -1ex \@minus -.2ex}%
                                   {2.3ex \@plus.2ex}%
                                   {\normalfont\normalsize\bf}}
\makeatother

\begin{document}

\begin{center}
{\Large \bf 
A cluster representation of the renewal Hawkes process
}
\end{center}
\footnotetext{
This research was supported by RIMS and by ISM.}
\begin{center}
Luis Iv\'an Hern\'andez Ru\'{\i}z\footnote{
Graduate School of Science, Kyoto University.}\footnote{
The research of this author was supported
by JSPS Open Partnership Joint Research Projects grant no. JPJSBP120209921.
}
\quad
and
\quad
Kouji Yano\footnotemark[2]\footnote{
Graduate School of Science, Osaka University.}\footnote{
The research of this author was supported by
JSPS KAKENHI grant no.'s JP19H01791, JP19K21834 and 21H01002.}
\end{center}

\begin{abstract}
A cluster representation for a Hawkes process with renewal immigration is obtained. The centre and satellite processes are indicated as a renewal process and generalized branching processes respectively. It is confirmed that the proposed construction indeed represents a cluster process and it is verified that it admits the desired intensity. Finally, the probability generating functional is computed for the stationary limit case.
\end{abstract}

\section{Introduction}
In Hawkes--Oakes \cite{HawkOak}, the linear Hawkes process was represented as a cluster process in which the centre process was given as a homogeneous Poisson process of \emph{immigrants} whose \emph{offspring} were given by branching satellite processes. Many generalizations of the linear Hawkes process have been studied in works such as Costa et. al. \cite{excinhibHawkes} in which a process with not only self-excitation, but also self-inhibition has been proposed. In Chen et. al. \cite{nonstationaryhawkes} a non-stationary version of the Hawkes process has been studied. Wheatley--Filimonov--Sornette \cite{Wheatley2016} generalized the linear Hawkes process by replacing the Poisson process with a renewal process and called it the \emph{renewal Hawkes process} abbreviated as RHP. Some methods for maximum likelihood estimation for the RHP were studied in Chen--Stindl \cite{ChenDirect} and furtherly refined in Chen--Stindl \cite{ChenImprove} to improve the speed of the computation. In this work, we obtain a cluster representation for the RHP and indicate explicitly what the centre and satellite processes are. We use a result of Westcott \cite{WestcottExistenceCluster} to show existence of the cluster process and verify that the construction indeed represents an RHP through its intensity. Finally, we find the limit process for the RHP at long times and compute its probability generating functional.

A \emph{simple point process} on $[0,\infty)$ is a sequence of nonnegative random variables $\{T_0,T_1,$ $T_2,\dots\}$ defined on a probability space $(\Omega,\cF,\bP)$ such that for all $n\ge0$, $T_{n+1}>T_n$ on $\cbra{T_n<\infty}$ and $T_{n+1}=\infty$ on $\cbra{T_n=\infty}$. We identify the point process with the associated counting process $N(t)=\sum_{i}1_{\cbra{T_i\le t}}$. Let $(\cF_t)$ be a filtration to which $N$ is adapted. An $\rbra{\cF_t}$-\emph{intensity} of $N$ is a nonnegative,  a.s. locally integrable process $\lambda_0(t)$ that is $(\cF_t)$-progressive, and such that,
\begin{align}
M(t)=N(t)-\int_0^t\lambda_0(s)\d s,
\end{align}
is an $(\cF_t)$-martingale. As a consequence, for any nonnegative process $C(t)$ which is \emph{predictable}, i.e. for all $t\ge0$, it is measurable with respect to the $\sigma$-field 
\begin{align}
\cP\rbra{\cF_t}=\sigma\rbra{(s,t]\times A; 0\le s\le t, A\in\cF_s}
\end{align}
it holds that,
\begin{align}
\label{integralIntensity}
\bE\sbra{\int_0^\infty C(s)N(\d s)}=\bE\sbra{\int_0^\infty C(s)\lambda_0(s)\d s}.
\end{align}
It is always possible (see e.g., \cite[Sec II. T12]{BremaudQueues}) to find a predictable version $\lambda$ of $\lambda_0$, in which case, it is essentially unique, i.e. $\lambda(t,\omega)=\lambda_0(t,\omega)$, $\bP(\d\omega)\d N(t,\omega)$-a.e. We define then the linear and renewal Hawkes processes through their intensity.

\begin{Def}A point process $N$ is called a \emph{linear Hawkes process} if $N$ admits a predictable  $\rbra{\cF_t}$-intensity given as
\begin{align}
\label{classicalHawkes}
\lambda(t)=\mu+\int_0^{t-}h(t-u)N(\d u),
\end{align}
where $\mu$ is a positive constant and $h$ is a nonnegative measurable function on $[0,\infty)$ satisfying $\int_0^\infty h(t)\d t<1$.
\end{Def}
Note that we have \eqref{integralIntensity} with 
\begin{align}
\lambda_0(t)=\mu+\int_0^t h(t-u)N(\d u),
\end{align}
because the integral with respect to the Lebesgue measure stays unaltered by adding one point at $t$. Thus $\lambda_0$ is an $(\cF_t)$-intensity of $N$, and $\lambda$ is the predictable version of $\lambda_0$.
 
It was shown in Hawkes--Oakes \cite{HawkOak} that the process with intensity \eqref{classicalHawkes} can be represented as a cluster process on $[0,\infty)$ with an homogeneous Poisson centre process of intensity $\mu$ and satellites given by generalized branching processes. These branching processes consist of inhomogeneous Poisson processes of characteristic intensity $h$ that start at each one of the previous points of the process up to time $t$.

The linear Hawkes process was generalized in Wheatley--Filimonov--Sornette \cite{Wheatley2016} by replacing the immigration Poisson process by a renewal process while keeping the structure of the offspring processes. In the same setting as above, we consider now a \emph{marked point process}, where the mark space is $\cbra{0,1}$ with its $\sigma$-field $2^{\cbra{0,1}}$ and the marks are random variables $D_i$, $i=1,2,\dots$, that take the values $D_i=0$ if the $i$-th point is an immigrant, and $D_i=1$ if it corresponds to an offspring.  The random variable $I(t)=\max\cbra{i;\;T_i\le t,D_i=0}$ represents the index of the last immigrant up to time $t$. We consider a filtration $(\cF_t)$ to which $N$ and $I$ are adapted. Additionally, consider a function $h$ satisfying the assumption:
\begin{enumerate}
\item[\textbf{(A)}]$h$ is a nonnegative measurable function on $[0,\infty)$ satisfying $\alpha:=\int_0^\infty h(t)\d t<1$.
\end{enumerate} 
We now give the definition of the RHP.
\begin{Def} A point process $N$ is called a \emph{renewal Hawkes process} if $N$ admits a predictable $(\cF_t)$-intensity given as
\begin{align}
\label{renewalHawkes}
\lambda(t)=\mu\rbra{t-T_{I(t-)}}+\int_0^{t-} h(t-u)N(\d u),
\end{align}
where $\mu$ is a nonnegative measurable function on $[0,\infty)$ such that 
\begin{align}
\label{hazardfunction}
\mu(t)=\frac{f(t)}{1-\int_0^t f(s)\d s},
\end{align}
for some probability density function $f$. The function $\mu$ is often called the \emph{hazard function}.
\end{Def}
The following lemma can help us get a better understanding of the role of the hazard function in the definition of the RHP. A proof is given in the appendix.
\begin{Lem}
\label{RenewalIntensity}
Let $\tau_1,\tau_2,\dots$ be i.i.d. random variables that have absolutely continuous distribution $F$ with density $f$. Define the \emph{renewals} $S_0=0$, $S_n:=\tau_1+\dots+\tau_n$, $n\ge 1$ and the associated counting process $N_R(t)=\sum_{i\ge0} 1_{\cbra{S_i\le t}}$, $t\ge0$. Let $(\cF_t)_{t\ge0}$ be a filtration to which $N_R$ is adapted. Then $N_R$ admits the predictable $(\cF_t)$-intensity 
\begin{align}
\mu(t-S_{N_R(t-)-1}),\quad t\ge0.
\end{align}
\end{Lem}
We see that the intensity for a process with i.i.d. interarrivals is given precisely by the hazard function shifted to the last arrival of such process. Notice that if we consider those interarrivals as exponentially distributed with parameter $1/\mu$ with $\mu>0$ a constant, then that makes the hazard function the constant $\mu$, and we recover the definition of the linear Hawkes process \eqref{classicalHawkes}. By adjusting the hazard function, we have control over the structure between immigrants.

The paper is organized as follows. A general framework for the theory of point processes and the definition of a cluster process are presented in Section \ref{SectionTheory}. In Section \ref{SectionConstruction} we propose a cluster representation for the RHP. In Section \ref{SectionExistence} we show that such a cluster process indeed exists. The verification that our representation admits the intensity of an RHP is discussed in Section \ref{SectionVerification}. Finally, Section \ref{SectionPgfl} and Section \ref{SectionStationary} are dedicated to the computation of the probability generating functional for the general RHP and the stationary limit process of the RHP respectively.

\section{Theoretical backgrounds for point processes}
In the following, the symbol of a measure $\nu$ on $[0,\infty)$ is used as well for its cumulative function $\nu(t)=\nu([0,t])$. Conversely, the symbol of a non-decreasing right-continuous function $\nu(t)$ on $[0,\infty)$ is used as well for its Stieltjes measure $\nu(\d t)$ such that $\nu(t)=\nu([0,t])$. For example, for a delta measure concentrated at a point $a\in\bR$, $\delta_a(\d t)$, the symbol $\delta_a(x)$ represents the function that is identically zero for $x<a$ and identically one for $x\ge a$.
\label{SectionTheory}
\subsection{The probability generating functional}
Let $\cX$ be a complete separable metric space and write $\cB(\cX)$ for the family of all Borel sets of $\cX$. Denote by $\Xi$ the class of measurable functions $z:\cX\rightarrow[0,1]$ such that $1-z$ vanishes outside some bounded set. A point process $N$ on $\cX$ is characterized in law (c.f. \cite[Corollary 9.2.IV]{DVJ2}) by a consistent family of finite dimensional distributions 
\begin{align}
\bP\rbra{\rbra{N(B_1),\dots,N(B_k)}\in \cdot},
\end{align}
for disjoint bounded sets $B_1,\dots,B_k$ in $\cB(\cX)$. The law of $N$ is then characterized by a \emph{probability generating functional} (p.g.fl.),
\begin{align}
\label{Defipgfl}
G[z]=\bE\sbra{\exp\rbra{\int_{\cX} \log z(t)N(\d t)}}=\bE\sbra{\prod_{t\in N(\cdot)}z(t)}\quad \text{for }z\in\Xi,
\end{align}
since for constants $\lambda_1,\dots,\lambda_k$ and disjoint bounded Borel sets $B_1,\dots,B_k$, by taking 
\begin{align}
z(\cdot)=\exp\rbra{-\sum_{i=1}^k\lambda_i 1_{B_i}(\cdot)},
\end{align}
we obtain
\begin{align}
G[z]=\bE\sbra{\exp\rbra{-\sum_{i=1}^k \lambda_i N(B_i)}},
\end{align}
which corresponds to the Laplace transform of $\rbra{N(B_1),\dots,N(B_k)}$. Notice that for the p.g.fl. to be well defined, one must make sure that the argument of the exponential function in \eqref{Defipgfl} is finite with probability one. This holds in general for point processes that are a.s. finite ($N(\cX)<\infty$ a.s.), but one could achieve this by restricting the functions $z$ to the class $\Xi$ (for details see Lemma 2 of Westcott \cite{Westcott1972}).

\subsection{Cluster processes}
Let $\cX$ and $\cY$ be complete separable metric spaces. Although cluster processes are defined in a general setting \cite[Section 6.3]{DVJ}, the case of interest is when $\cX=\cY=\bR$, so there is no loss of understanding if the reader assumes  $\cX=\cY=\bR$. Let $\cN_{\cX}^\sharp$ stand for the set of counting measures $\nu$ on $\cX$ which are \emph{locally finite} in the sense that $\nu(B)<\infty$ for all bounded set $B\in\cB(\cX)$. We write $\cB(\cN_{\cX}^\sharp)$ to denote the $\sigma$-field generated by all the maps $\cN_{\cX}^\sharp\ni\nu\mapsto\nu(B)$ for $B$ in $\cB(\cX)$, and $\cP(\cN_{\cX}^\sharp)$ to denote the space of probability measures on $\cN_{\cX}^\sharp$. The \emph{convolution} of $\Pi$ and $\Pi^\prime\in\cP(\cN_{\cX}^\sharp)$ is defined as
\begin{align}
\label{ConvLaw}
\rbra{\Pi*\Pi^\prime}(U)=\mathrel{\mathop{\int}\limits^{}_{\cN_{\cX}^\sharp\times\cN_{\cX}^\sharp}}1_{\rbra{\nu+\nu^\prime\in U}}\Pi(\d\nu)\Pi^\prime(\d\nu^\prime)\quad\text{for all }U\in\cB(\cN_{\cX}^\sharp).
\end{align}
If we have two independent point processes $N(\cdot)$ and $N^\prime(\cdot)$ on $\cX$, we can then use \eqref{ConvLaw} to write the law of their sum. Let $\Pi(\cdot)=\bP\rbra{N\in\cdot}$ and $\Pi^\prime(\cdot)=\bP\rbra{N^\prime\in\cdot}$, then
\begin{align}
\Pi*\Pi^\prime=\bP(N+N^\prime\in\cdot).
\end{align}
We can also write an expression for the p.g.fl. If we denote
\begin{align}
G_\Pi[z]=\bE\sbra{\exp\int\log z(t)N(\d t)},\quad\text{and }\;  G_{\Pi^\prime}[z]=\bE\sbra{\exp\int\log z(t)N^\prime(\d t)},
\end{align}
where the expectation is taken w.r.t. $\Pi$ and $\Pi^\prime$ respectively, then the p.g.fl. of $N+N^\prime$ is given as
\begin{align}
G_{\Pi*\Pi^\prime}[z]=&\bE\sbra{\exp\int\log z(t) (N+N^\prime)(\d t)}
\label{}\\
=&G_{\Pi}[z]G_{\Pi^\prime}[z].
\label{PgflSum}
\end{align}

\begin{Def}[\textbf{\cite[Ch.6 p.165]{DVJ}}]
A (symbolic) \emph{measurable family of point processes} on $\cX$ is a family $\{N(\cdot \mid y):y\in \cY \}$ where for all $y\in\cY$, $N(\cdot \mid y)$ is a point process on $\cX$, and for all $U\in\cB(\cN_{\cX}^\sharp)$ the function 
\begin{align}
y\longmapsto \bP(N(\cdot\mid y)\in U)
\label{}
\end{align}
is $\cB(\cY)$-measurable.
\end{Def}

The construction of a cluster process involves two components \cite{EarthquakesVJ}: a point process $N_c$ of cluster centres whose realization consists of the points $\cbra{y_i}_{i\ge0} \subset \cY$, and a family of point processes on $\cX$, namely $\cbra{N_s(\cdot \mid y); y\in\cY}$, whose superposition constitute the observed process. We formalize this idea through the convolution in $\cP(\cN_{\cX}^\sharp)$.
\begin{Def}
Let $N_c$ be a point process on $\cY$ and $\left\{N_s(\cdot \mid y): \right.$ $\left. y\in \cY\right\}$ a measurable family of point processes on $\cX$. (The family $\left\{N_s(\cdot \mid y): \right.$ $\left. y\in \cY\right\}$ is considered to be mutually independent and to be independent of $N_c$.) Then, the \emph{independent cluster process} on $\cX$, with \emph{centre process} $N_c$ and \emph{satellite processes} $\{N_s(\cdot \mid y): y\in\cY\}$, which we denote by
\begin{align}
N(\cdot)=\int_\cY N_s(\cdot \mid y)N_c(dy) = \sum_{y \in N_c(\cdot)} N_s(\cdot \mid y),
\label{ClusterDefinition}
\end{align}
is defined in law as
\begin{align}
\label{ClusterLaw}
\bP(N\in U)=\int_{N^\sharp_{\cY}}\bP\rbra{N_s(\cdot\mid\mu)\in U}\bP\rbra{N_c\in\d\mu},
\end{align}
where $\bP\rbra{N_s(\cdot\mid\mu)\in U}$ for $\mu(\cdot)=\sum_{i}\delta_{y_i}(\cdot)\in\cN^\sharp_{\cY}$ is defined as the infinite convolution
\begin{align}
\bP\rbra{N_s(\cdot\mid\mu)\in U}=\rbra{\Pi_{y_1}*\Pi_{y_2}*\cdots}(U)\quad\text{for } U\in\cB(\cN^\sharp_{\cX}),
\end{align}
with $\Pi_y(U)=\bP\rbra{N_s(\cdot\mid y)\in U}$ for $U\in\cB(\cN^\sharp_{\cX})$ and $y\in\cY$.
\end{Def}
We now give an expression for the p.g.fl. of the independent cluster process in the following Theorem (see for example, \cite[Equation (3) of discussion by Moyal on page 37]{NeymanScottCosmology}). For completion of this paper, a proof is included in the appendix.

\begin{Thm}
\label{TheoremPgfl}
Let $N$ be an independent cluster process with centre process $N_c$ and satellite processes $\left\{N_s(\cdot \mid y): \right.$ $\left. y\in \cY\right\}$. Let $G_c[z]$ denote the p.g.fl. of the centre process and $G_s[z\mid y]$ the p.g.fl. of $N_s(\cdot \mid y)$. Then $G[z]$, the p.g.fl. of $N(\cdot)$, is given by
\begin{align}
\label{newpgfl}
G[z]=G_c\sbra{G_s[z\mid \cdot]}\quad\text{for any }z\in\Xi.
\end{align}
\end{Thm}

A necessary and sufficient condition for the superposition \eqref{ClusterDefinition} to define a point process on $\cX$ \cite[Ch. 6 Equation (6.3.1)]{DVJ}, in which case we say that the independent cluster process exists, is that for every bounded set $B\in \cB(\cX)$,
\begin{align}
N(B)=\int_\cY N_s(B \mid y)N_c(dy) = \sum_{y \in N_c(\cdot)} N_s(B \mid y)<\infty\quad\text{a.s.}
\label{ClusterCondition}
\end{align}
Equivalent conditions for verifying \eqref{ClusterCondition} were presented in Westcott \cite{WestcottExistenceCluster}, namely the following theorem and its corollary.
\begin{Thm}[\textbf{\cite[Theorem 3]{WestcottExistenceCluster}}]
\label{WestcottThm}
Under the assumptions that the clusters are identically distributed, i.e.
\begin{align}
N_s(\cdot\mid y)\dist N_s(\cdot-y\mid 0)\quad \text{for all } y\in\cY,
\end{align} 
and that the clusters are almost surely finite, namely
\begin{align}
N_s(\cX\mid 0)<\infty\quad\text{a.s.},
\end{align}
the independent cluster process $N$ exists if and only if for every bounded set $B\in\cB(\cX)$,
\begin{align}
\label{FiniteClusterCondition}
\int_{\cY}\bP\rbra{N_s(B\mid y)>0}N_c(\d y)<\infty \quad a.s.
\end{align}
\end{Thm}

We give its proof in the appendix for completeness of this paper.

\begin{Cor}[\textbf{\cite[Corollary 3.3]{WestcottExistenceCluster}}]
\label{WestcottCor}
Let $\cX,\cY=\bR$ and assume the following conditions:
\begin{enumerate}
\item[\emph{(i)}]$\sup_{t}\bE\sbra{N_c(I-t)}<\infty$ for all bounded interval $I$.
\item[\emph{(ii)}]$N_s(\cdot\mid t)\dist N_s(\cdot-t\mid 0)$ for all $t\in\bR$.
\item[\emph{(iii)}]$\bE\sbra{N_s(\bR\mid 0)}<\infty$.
\end{enumerate}
Then \eqref{FiniteClusterCondition} is satisfied.
\end{Cor}
We give its proof in the appendix for completeness of this paper.

\section{Construction of the RHP}
Understanding the cluster structure of the RHP plays an important role in proving limit theorems, as was done in Hern\'andez \cite{IvanLimits2023} by allowing to treat the immigration process with the tools of renewal theory, or in \cite{BordenaveLarge} where a Large Deviation Principle was derived from the cluster representation of the linear Hawkes process. The cluster representation has also been used for perfect simulation of the linear Hawkes process \cite{ChenSim}. 

By changing the form of the immigration process, we have the liberty to adjust the spacing between clusters, which gives more freedom for application purposes. While the family trees formed by the offspring of the immigrants have the same structure as those presented by Hawkes--Oakes \cite{HawkOak}, in that work, they miss to provide a mathematical expression of the process as a superposition of the clusters generated by each immigrant. In this section, we obtain such a formula and use it for the computation of relevant quantities in subsequent sections. 

Before proceeding, we make a remark on the notation used.

\label{SectionConstruction}
If $f$ and $g$ are both functions, we will denote their convolution as 
\begin{align}
f*g(t)=\int_0^tf(t-s)g(s)\d s,
\end{align}
whereas, if $F$ is a measure on $[0,\infty)$ and $g$ is a function, the convention that $F*g$ is a function is used, and we write
\begin{align}
F*g(t)=\int_0^t g(t-s)F(\d s).
\end{align}
We sometimes identify the measure $F(\d s)$ with its cumulative distribution function $F(t)=\int_0^t F(\d s)$. If $F$ and $G$ are both measures, we will denote their convolution as
\begin{align}
F*G(t)=\int_0^t F(t-s)G(\d s)=\int_0^t G(t-s)F(\d s).
\end{align} 

\subsection{The centre process}
In the RHP, immigration is given by a \emph{renewal process}, which naturally we take as our centre process. Let $\tau,\tau_1,\tau_2,\dots,$ be positive i.i.d. random variables whose probability distribution function
\begin{align}
\label{renewalDistribution}
F(t):=\bP(\tau\le t),
\end{align} 
satisfies the assumption:
\begin{enumerate}
\item[\textbf{(B)}]$\tau$ has finite mean $m^{-1}:=\bE\sbra{\tau}$ and its distribution function $F$ has a density $f$, i.e. $F(t)=\int_0^t f(s)\d s$.
\end{enumerate} 
The epochs of the renewal process correspond to the partial sums $S_0=0$, $S_n=\tau_1+\dots+\tau_n$, to which we associate the counting process $N_R(t)=\sum_{i\ge 0} 1_{\cbra{S_i\le t}}$. For each $n\ge0$ and $x\ge0$, the distribution of $S_n$ is $\bP(S_n\le x)=F^{*n}(x)$, where

\begin{align}
F^{*0}(x)=\delta_0(x),\quad F^{*(n+1)}(x)=\int_0^xF^{*n}(x-y)F(\d y).
\end{align}

The mean number of events up to time $t$ is given by the \emph{renewal function},
\begin{align}
\label{RenewalFunction}
\Phi(t):=\bE\sbra{N_R(t)}=\sum_{n\ge0}F^{*n}(t).
\end{align}
Note that $\Phi$ is an increasing function. Since the distribution function $F$ has a density $f$, the \emph{renewal measure} $\Phi(\d t)$ is absolutely continuous w.r.t. the Lebesgue measure with density $\varphi(t)=\sum_{n\ge1}f^{*n}(t)$ for all $t\ge0$ \cite[Proposition V.2.7]{Asmussen2003}.

Denote by $(\cF^R_t)_{t\ge0}$ the augmentation of the natural filtration of the renewal process $\sigma\rbra{N_R(s);0\le s\le t}$, and take the hazard function $\mu(t)$ as in \eqref{hazardfunction}. We use Lemma \ref{RenewalIntensity} to see that
\begin{align}
\lambda_R(t):=\mu(t-S_{N_R(t-)-1}),
\end{align}
is the predictable $(\cF^R_t)$-intensity of $N_R$.

\subsection{The satellite processes}
One feature of the RHP is that in addition to immigrants being able to generate offspring, these offspring themselves can generate further offspring. Thus, offspring points could be described as higher-level centre processes. In order to represent this structure, we use a notation similar to Neyman--Scott \cite{NeymanScottCosmology}. The renewal process will be named a zero-th order centre process $N_c^{(0)}(\cdot):=N_R(\cdot)$. 

Let $\cbra{N_s^{\rbra{n}}(\cdot\mid t);\;t\in[0,\infty)}_{n\ge1}$ be a sequence of measurable families of point processes which is considered to be i.i.d. and to be independent of $N_c^{(0)}$, such that $N_s^{(n)}(\cdot\mid t)$ has the same law as $N_s(\cdot\mid t)$ whose p.g.fl. is given by
\begin{align}
\label{Defpgfl}
\bE\sbra{\prod_{x\in N_s(\cdot\mid t)}z(x)}=\exp\rbra{\int_0^\infty\rbra{z(x+t)-1}h(x)\d x},
\end{align}
with $h$ satisfying \textbf{(A)}.
If we take $z(x)=e^{-\lambda}$ for $\lambda>0$, we have
\begin{align}
\bE\sbra{e^{-\lambda N_s\rbra{[0,\infty)\mid t}}}=\exp\rbra{\rbra{e^{-\lambda}-1}\int_0^\infty h(x)\d x}=e^{-\alpha\;\rbra{1-e^{-\lambda}}},
\end{align}
which shows that $N_s([0,\infty)\mid t)\dist\text{Poi}(\alpha)$. More generally, take disjoint intervals $(a_1,b_1],$ $\dots,$ $(a_k,b_k]$. If we substitute the function
\begin{align}
z(x)=\left\{\begin{array}{ll} e^{-\lambda_i}&\rbra{x\in(a_i,b_i]},\\1&\rbra{x\notin \cup_{i=1}^k(a_i,b_i]}, \end{array}\right.
\end{align}
into \eqref{Defpgfl}, we obtain
\begin{align}
\bE\sbra{e^{-\lambda_1 N_s((a_1,b_1]\mid t)}\cdots e^{-\lambda_k N_s((a_k,b_k]\mid t)}}=\prod_{i=1}^k\exp\rbra{(e^{-\lambda_i}-1)\int_{a_i}^{b_i} h(x-t)\d x}.
\end{align}
This shows that $\cbra{N_s((a_i,b_i]\mid t)}_{i=1}^k$ is mutually independent and 
\begin{align}
\label{PoissonComponent}
N_s((a_i,b_i]\mid t)\dist\text{Poi}\rbra{\int_{a_i}^{b_i}h(x-t)\d x},
\end{align}
or in other words, $N_s(\cdot\mid t)$ is an inhomogeneous Poisson process with intensity $h(x-t)\d x$, where we understand $h(x)=0$ for $x<0$. In particular, we see that $\cbra{N_s(\cdot\mid t): t\ge0}$ is a measurable family of point processes.

Given that there is a centre at $t_0\ge 0$, we construct higher-level centre processes $N_c^{(n)}$, for $n\ge1$, from a superposition of the processes $N_s^{(n)}$ with the following recursive structure:
\begin{align}
\label{}
N_c^{(0)}(\cdot\mid t_0):=\delta_{t_0},\quad N_c^{\rbra{n+1}}(\cdot\mid t_0)=\sum_{t\in N_c^{\rbra{n}}(\cdot\mid t_0)}N_s^{\rbra{n+1}}(\cdot\mid t),
\end{align}
where $N_c^{(0)}(\cdot\mid t_0)$ is the original immigrant at $t_0$ and $N_c^{\rbra{n}}(\cdot\mid t_0)$ represents its $n$-th generation offspring. We define as well some processes of interest, namely, the total number of $n$-th generation descendants,
\begin{align}
N_c^{(n)}(\cdot)=\sum_{t_0\in N_R(\cdot)} N_c^{(n)}(\cdot\mid t_0),
\end{align}
and the complete offspring of the imimigrant at $t_0$ (including the immigrant),
\begin{align}
\label{ClusterComponent}
N_c(\cdot\mid t_0)=\sum_{n\ge0}N_c^{(n)}(\cdot\mid t_0).
\end{align}
 We take the processes defined as in $\eqref{ClusterComponent}$ as the satellite processes of our construction for a centre located at $t_0$. Finally, the RHP is given by the superposition:
\begin{align}
\label{RHPCluster}
N(\cdot)=\int_0^\infty N_c(\cdot\mid t)N_R(\d t)=\sum_{t_0\in N_R(\cdot)}\sum_{n\ge0}N_c^{(n)}(\cdot\mid t_0).
\end{align}
Note that $\eqref{RHPCluster}$ can also be written as
\begin{align}
\label{RHPAlternative}
N(\cdot)=\sum_{n\ge0}\sum_{t_0\in N_R(\cdot)}N_c^{(n)}(\cdot\mid t_0)=N_R(\cdot)+\sum_{n\ge1}N_c^{(n)}(\cdot).
\end{align}

\section{Validity of the construction}
\label{SectionExistence}
We are concerned with whether our construction of the RHP in fact represents a valid cluster process, since in the original work by Wheatley--Filimonov--Sornette \cite{Wheatley2016} no result of existence was provided. In the following, we show that with the assumptions made for its construction, \eqref{RHPCluster} is a valid definition.

\begin{Thm}
\label{RHPExistenceGeneral}
Let $N_g$ be a point process on $[0,\infty)$ such that
\begin{align}
\label{UnifBoundCenter}
\sup_{t}\bE\sbra{N_g(I-t)}<\infty,\quad \text{for all bounded interval } I\subset[0,\infty),
\end{align}
and that is independent from $\cbra{N_s^{\rbra{n}}(\cdot\mid t);\;t\in[0,\infty)}_{n\ge1}$. If $h$ is a function satisfying \emph{\textbf{(A)}}, then the cluster process defined as
\begin{align}
\sum_{t_0\in N_g(\cdot)}\sum_{n\ge0}N_c^{(n)}(\cdot\mid t_0),
\end{align}
with immigrant process $N_g$, exists and has a.s. finite clusters.
\end{Thm}
\Proof{
Using Corollary \ref{WestcottCor} , the existence of the cluster process is proved if we can show that conditions (i)-(iii) hold.

Claim (i) is satisfied by assumption.

Claim (ii) follows from the construction of the satellites as superposition of inhomogeneous Poisson processes with p.g.fl. \eqref{Defpgfl} that originate at previous points of the process and the observation \eqref{PoissonComponent}.

Finally, to prove claim (iii), let $t_0\ge0$ and consider $N_c(\cdot\mid t_0)$ a cluster with centre at $t_0$. Let us call $Z_n:=N_c^{(n)}([0,\infty)\mid t_0)$ for $n\ge0$ and $\bar{Z}:=N_c([0,\infty)\mid t_0)$ so that $Z_n$ and $\bar{Z}$ represent respectively the number of $n$-th generation points and the total number of points in the cluster. These random variables satisfy:
\begin{align*}
Z_n=\sum_{k=0}^{Z_{n-1}}X_{n,k},
\label{}
\end{align*}
$N_c^{(n)}(\cdot\mid t_0)=\{Y_1^{(n)},\dots,Y_{Z_n}^{(n)}\}$ for $n\ge0$, and $X_{n,k}:=N_s^{(n)}([0,\infty)\mid Y_k^{(n-1)})$. This shows that the variables $Z_n$ form a Galton--Watson process with offspring density $h$, and that the number of points per generation follows a Poisson distribution of parameter $\alpha$. Then $\bar{Z}$ is given by the total size of the Galton--Watson process. From the Galton--Watson theory \cite[Theorem 6.1, p.7]{Harris1963} we know that $\bP\rbra{\bar{Z}<\infty}=1$ if $\alpha<1$, and the p.g.f. of the cluster size $\pi(u)=\bE\sbra{u^{\bar{Z}}}$, $0<u<1$   \cite[Section 13.2, p.32]{Harris1963} satisfies,
\begin{align*}
\pi(u)=u\exp\cbra{\alpha\sbra{\pi(u)-1}},
\end{align*} 
from which we can conclude that $\pi^\prime(1-)=\bE\sbra{\;\bar{Z}\;}=\frac{1}{1-\alpha}$. We then see that the three assumptions needed are satisfied, which concludes the proof.
}

We now see that our construction \eqref{RHPCluster} is a cluster process.

\begin{Cor}
\label{RHPExistence}
If $h$ is a function satisfying \emph{\textbf{(A)}}, and $N_R(\cdot)$ is a renewal process satisfying \emph{\textbf{(B)}}, then the cluster process defined as in \eqref{RHPCluster} exists and has a.s. finite clusters.
\end{Cor}
\Proof{
Using Theorem \ref{RHPExistenceGeneral} it is enough to note that $N_R$ satisfies \eqref{UnifBoundCenter}, which follows from the bound,
\begin{align}
\label{RenewalSubadditive}
\bE\sbra{N_R(t+a)-N_R(t)}\le\Phi(a),\quad\text{for all }a>0,
\end{align}
(see for example \cite[Sec. V, Theorem 2.4, p.146]{Asmussen2003}). Let $I$ be any bounded interval in $[0,\infty)$ and write $\absol{I}$ for its Lebesgue measure. Then, from \eqref{RenewalSubadditive} we can see that $\bE\sbra{N_R(I)}\le\Phi(\absol{I})$ and we can conclude that,
\begin{align}
\supre{}{t\ge0}\bE\sbra{N_R(I-t)}\le\supre{}{t\ge0}\Phi\rbra{\absol{I-t}}=\Phi(\absol{I})<\infty.
\end{align}
This completes the argument.

The cluster structure of the linear Hawkes process has been thoroughly studied. For example, Reynaud-Bouret--Roy \cite{ReynaudNonAs} found tail estimates for the length of a cluster and for the extinction time of a cluster, and Daw \cite{DawCluster} computed the distribution for the length of a cluster. Because the branching structure of the processes treated in Theorem \ref{RHPExistenceGeneral} is the same as that of the linear Hawkes process, then results on the properties of the clusters, like the ones mentioned, generalize directly to these new processes.
}

\section{Verification of the intensity}
\label{SectionVerification}
In this section we show that the process constructed in the previous section admits the desired intensity of an RHP. We point out that in Hawkes--Oakes \cite[Lemma 1]{HawkOak}, an informal argument was given for the case of the linear Hawkes process by appealing to the structure of an age-dependent pure birth process. Here we give a formal argument appealing to the uniqueness of predictable intensities.
\begin{Thm}
\label{RHPValidity}
Assume that $h$ is a function satisfying \emph{\textbf{(A)}}, and that $N_R(\cdot)$ is a renewal process that satisfies assumption \emph{\textbf{(B)}}, then the cluster process defined as in \eqref{RHPCluster} has the predictable intensity \eqref{renewalHawkes}.
\end{Thm}
\Proof{
Define the random variables
\begin{align}
D_i=\left\{\begin{array}{ll} 0\quad&\rbra{\text{if }T_i\in N_R(\cdot)},\\1&\rbra{\text{otherwise}}, \end{array}\right.
\end{align}
and the function $I(t)=\max\cbra{i;\;T_i\le t,D_i=0}$ as before. We construct the filtration $(\cF_t)_{t\ge0}$ by the augmentation of the natural filtration $(\cF^0_t)_{t\ge0}$ defined as
\begin{align}
\cF^0_t=\sigma\rbra{N_c^{(n)}\rbra{(a,b]}; \;0\le a\le b\le t, \;n=0,1,2,\dots},\quad t\ge0.
\end{align}

Notice that in the definition of the intensity \eqref{renewalHawkes}, the term $\mu\rbra{t-T_{I(t-)}}$ is only affected by the terms that come from the renewal process. This is because the points $T_{I(t-)}$ all correspond to immigrants. Since $N_R(t)=\min\cbra{i:S_i\le t}$, let us then denote $N_R(\cdot)=\cbra{T_i: D_i=0}=\cbra{S_{1},S_{2},\dots}$, and notice then that 
\begin{align}
\label{renewalInt1}
\mu(t-T_{I(t-)})=\mu(t-S_{N_R(t-)-1}),\quad t\ge0.
\end{align}
Consider now an arbitrary $(\cF_t)$-predictable process $C(u)=1_A1_{(r,t]}(u)$ for $A\in\cF_r$ and $0\le r\le t$. Then
\begin{align}
\bE\sbra{\int C(u)N(\d u)}=\bE\sbra{1_A\int_r^tN_R(\d u)}+\bE\sbra{1_A\int_r^t\sum_{n\ge0}N^{(n+1)}_c(\d u)}.
\end{align}
Since $\lambda_R(t)$ is an $(\cF_t)$-intensity of $N_R$, the first term on the R.H.S. equals
\begin{align}
\bE\sbra{1_A\int_r^tN_R(\d u)}=\bE\sbra{1_A\int_r^t\mu(u-S_{N_R(u-)})\d u}=\bE\sbra{1_A\int_r^t\mu(u-T_{I(u-)})\d u}.
\end{align}
We also have,
\begin{align}
\bE\sbra{1_A\int_r^t\sum_{n\ge0}N_c^{(n+1)}(\d u)}=&\bE\sbra{1_A\int_r^t\sum_{n\ge0}\sum_{x\in N_c^{(n)}(\cdot)}N_s^{(n+1)}(\d u\mid x)}
\label{}\\
=&\sum_{n\ge0}\bE\sbra{1_A\sum_{x\in N_c^{(n)}(\cdot)}N_s^{(n+1)}((r,t]\mid x)}.
\end{align}
Because of the independence property of a Poisson point process, we have independence of $\cbra{N_s^{\rbra{n+1}}((r,t]\mid x);\;t>r,x\ge0}$ from $\cF_r\vee\sigma\rbra{N_c^{(n)}}$, and by \eqref{PoissonComponent} we obtain
\begin{align}
\sum_{n\ge0}\bE\sbra{1_A\sum_{x\in N_c^{(n)}(\cdot)}N_s^{(n+1)}((r,t]\mid x)}=&\sum_{n\ge0}\bE\sbra{1_A\sum_{x\in N_c^{(n)}(\cdot)}\int_r^th(u-x)\d u}
\label{}\\
=&\bE\sbra{1_A\int_0^\infty\rbra{\int_r^th(u-x)\d u}\sum_{n\ge0}N_c^{(n)}(\d x)}
\label{}\\
=&\bE\sbra{1_A\int_0^\infty\rbra{\int_r^th(u-x)\d u}N(\d x)}
\label{}\\
=&\bE\sbra{1_A\int_r^t\rbra{\int_0^uh(u-x)N(\d x)}\d u}.
\end{align}
From the chain of equalities we obtain the identity
\begin{align}
\bE\sbra{\int C(u)N(\d u)}=\bE\sbra{\int C(u)\cbra{\mu(u-T_{I(u-)})+\int_0^{u-}h(u-x)N(\d x)}\d u},
\label{VerificationIntensity}
\end{align}
which may be extended to all $(\cF_t)$-predictable processes $C$. Thus, we can say that $N$ admits the predictable $(\cF_t)$-intensity
\begin{align}
\lambda(t):=\mu(t-T_{I(t-)})+\int_0^{t-}h(t-x)N(\d x),\quad t\ge0.
\end{align} 
}

\textbf{Remark.} In order to appeal to the results of Westcott \cite{WestcottExistenceCluster} in the proof of Theorem \ref{RHPExistenceGeneral}, we have made use of assumption \textbf{(A)}, in particular, the condition that $\int_0^\infty h(t)\d t < 1$. In the case of the linear Hawkes process, this condition can be relaxed to prove non-explosiveness by asking only local integrability of the function $h$, i.e. that for any $T>0$, the integral $\int_0^T h(t)\d t<\infty$, and use a result of Massoulié \cite{Massoulie1998} to conclude. Nevertheless, this strategy cannot be directly applied to the RHP without making further regularity assumptions on $\mu(x)$, since it becomes necessary to bound $\absol{\mu(t-T_{I(t)})-\mu(t-T^\prime_{I^\prime(t)})}$ where $T$ and $T^\prime$ correspond to different immigration processes with the same inter-arrival distribution. The fact that we can obtain non explosiveness properties from the cluster construction emphasizes the importance of having an explicit cluster representation.

\section{Probability generating functional for the RHP}
\label{SectionPgfl}
In this section we investigate the p.g.fl. of the RHP. The difficulty of finding the complete p.g.fl. lies on the fact that the renewal process is not a finite point process. We begin with the p.g.fl. of the renewal process. 

Denote $p_n(T):=\bP\rbra{N_R((0,T])=n}$ for $T\ge0$ and any nonnegative integer $n$, and for any $z\in\Xi$ and $T\ge0$ define
\begin{align}
z^T(t):=\left\{\begin{array}{lc}z(t)&(t\le T),\\1&(t>T).\end{array}\right.
\end{align}
Since $\cbra{z^T}_{T\ge0}$ is decreasing in $T$ for $z\in\Xi$, and all $z^T$ are dominated by the constant $1$, it is a consequence of the DCT that
\begin{align}
G_R[z]=\lim_{T\rightarrow\infty}G_R[z^T],
\end{align}
where
\begin{align}
&G_R[z^T]=\sum_{n=0}^\infty\bE\sbra{\prod_{x\in N_R(\cdot)}z^T(x)\;; N_R((0,T])=n}
\label{}\\
=&p_0(T)+\sum_{n=1}^\infty\bE\sbra{z(S_1)\cdots z(S_n)\;; N_R((0,T])=n}
\label{}\\
=&p_0(T)+\sum_{n=1}^\infty p_n(T)\int_0^T\int_0^{T-s_1}\cdots\int_0^{T-s_{n-1}}z(s_1)\cdots z(s_n)f(s_1)\cdots f(s_n-s_{n-1})\d s_n\cdots\d s_1,
\end{align}
but this expression cannot be furtherly simplified in general. 

Now we focus on the p.g.fl. for the satellite processes. A cluster whose centre is located at $t_0$, for some $t_0\ge0$, is formed by the immigrant that originated it together with all the generations of its offspring. The following relation has been already established in Hawkes--Oakes \cite{HawkOak}, and we provide a proof using our construction.
\begin{Thm}
Let $t_0\ge0$. The p.g.fl. for a cluster starting at $t_0$, namely 
\begin{align}
G_c[z\mid t_0]=\bE\sbra{\exp\int_0^\infty\log z(t)N_c(\d t\mid t_0)},
\end{align}
satisfies the functional relation
\begin{align}
\label{pgflComponent}
G_c\sbra{z\mid t_0}=z(t_0)\exp\cbra{\int_0^\infty \rbra{G_c[z(x+\cdot)\mid t_0]-1}h(x)\d x}.
\end{align}
\end{Thm}
\Proof{
Let $t_0\ge0$. We want the p.g.fl. for the cluster 
\begin{align}
N_c(\cdot\mid t_0)=\sum_{n\ge0}N_c^{(n)}(\cdot\mid t_0),
\end{align} 
Let us call $G_c^{(n)}$ the p.g.fl. of the cluster up to generation $n$, namely 
\begin{align}
G_c^{(n)}[z\mid t_0]=\bE\sbra{\exp\int_0^\infty\log z(t)\sum_{i=0}^n N_c^{(i)}(\d t\mid t_0)}. 
\end{align}
Denote $\cF_0:=\cbra{\Omega,\emptyset}$ and $\cF_n:=\sigma\rbra{N_c^{(1)}(\cdot\mid t_0),\dots,N_c^{(n)}(\cdot\mid t_0)}$.
We have,
\begin{align}
G_c^{(n+1)}[z\mid t_0]=&\bE\sbra{\exp\int\log z(x)\sum_{i=0}^{n+1}N_c^{(i)}(\d x\mid t_0)}
\label{}\\
=&\bE\sbra{\exp\int\log z(x)\rbra{\delta_{t_0}(\d x)+\sum_{i=1}^{n+1}N_c^{(i)}(\d x\mid t_0)}}
\label{}\\
=&z(t_0)\bE\sbra{\exp\cbra{\int\log z(x)\sum_{i=0}^{n}N_c^{(i+1)}(\d x\mid t_0)}}.
\end{align}
Note that the above expectation can be written as,
\begin{align}
&\bE\sbra{\exp\cbra{\int\log z(x)\sum_{i=0}^{n}N_c^{(i+1)}(\d x\mid t_0)}}
\label{}\\
=&\bE\sbra{\exp\cbra{\sum_{i=0}^{n}\sum_{y\in N_c^{(i)}(\cdot\mid t_0)}\int\log z(x)N_s^{(i+1)}(\d x\mid y)}}
\end{align}
Using that the processes $N_s^{(i+1)}(\cdot\mid y)$ are an i.i.d. family and independent from $N_c^{(i)}(\cdot\mid t_0)$, we can rewrite the above expectation as,

\begin{align}
&\bE\sbra{\exp\cbra{\sum_{i=0}^{n}\sum_{y\in N_c^{(i)}(\cdot\mid t_0)}\int\log z(x)N_s^{(i+1)}(\d x\mid y)}}
\label{}\\
=&\bE\sbra{\prod_{i=0}^{n}\prod_{y\in N_c^{(i)}(\cdot\mid t_0)}\bE\sbra{\exp\cbra{\int\log z(x)N_s^{(i+1)}(\d x\mid y)}\bigg| \cF_i}}
\label{TowerAnverse}\\
=&\bE\sbra{\prod_{i=0}^{n}\prod_{y\in N_c^{(i)}(\cdot\mid t_0)}\bE\sbra{\exp\cbra{\int\log z(x)N_s^{(i+1)}(\d x\mid y)}}},
\label{}
\end{align}
where the tower property was recursively applied to obtain \eqref{TowerAnverse}. Now we use that for any $i\ge0$ and $t\ge0$, $N_s^{(i)}(\cdot\mid t)\dist N_s^{(i)}(\cdot-t\mid 0)$. Furthermore, since the 
\begin{align}
\{N_s^{(i)}(\cdot-t\mid 0)\}_{i\ge0},
\end{align} 
form an i.i.d. sequence of point processes, we can exploit the independence of $N_s^{(n+1)}(\cdot-t\mid 0)$ from the $\{N_c^{(i)}(\cdot-t\mid 0)\}_{i=0}^{n}$ and write
\begin{align}
&\bE\sbra{\prod_{i=0}^{n}\prod_{y\in N_c^{(i)}(\cdot\mid t_0)}\bE\sbra{\exp\cbra{\int\log z(x)N_s^{(i+1)}(\d x\mid y)}}}.
\label{}\\
=&\bE\sbra{\prod_{i=0}^{n}\prod_{y\in N_c^{(i)}(\cdot\mid t_0)}\bE\sbra{\exp\cbra{\int\log z(x+y)N_s^{(n+1)}(\d x\mid 0)}}}
\label{}\\
=&\bE\sbra{\prod_{i=0}^{n}\prod_{y\in N_c^{(i)}(\cdot\mid t_0)}\bE\sbra{\exp\cbra{\int\log z(x+y)N_s^{(n+1)}(\d x\mid 0)}\bigg| \cF_i}}
\label{TowerReverse}\\
=&\bE\sbra{\prod_{i=0}^{n}\prod_{y\in N_c^{(i)}(\cdot\mid t_0)}\exp\cbra{\int\log z(x+y)N_s^{(n+1)}(\d x\mid 0)}}
\label{}\\
=&\bE\sbra{\exp\cbra{\sum_{i=0}^{n}\sum_{y\in N_c^{(i)}(\cdot\mid t_0)}\int\log z(x+y)N_s^{(n+1)}(\d x\mid 0)}},
\end{align}
where \eqref{TowerReverse} was obtained by once again recursively applying the tower property, but in reverse order compared to \eqref{TowerAnverse}. Writing the summation above as an integral w.r.t. the counting measure $N_c^{(i)}(\cdot\mid t_0)$ and using Fubini's Theorem we get,
\begin{align}
&\bE\sbra{\exp\cbra{\sum_{i=0}^{n}\sum_{y\in N_c^{(i)}(\cdot\mid t_0)}\int\log z(x+y)N_s^{(n+1)}(\d x\mid 0)}}
\label{}\\
=&\bE\sbra{\exp\cbra{\int\int\log z(x+y)\sum_{i=0}^{n}N_c^{(i)}(\d y\mid t_0)N_s^{(n+1)}(\d x\mid 0)}}
\label{}\\
=&\bE\sbra{\prod_{x\in N_s^{(n+1)}(\cdot\mid 0)}\bE\sbra{\exp\cbra{\int\log z(x+y)\sum_{i=0}^{n}N_c^{(i)}(\d y\mid t_0)}}}.
\label{}
\end{align}
This last step follows from the independence of $N_s^{(n+1)}(\cdot\mid 0)$ from the $\{N_c^{(i)}(\cdot\mid t_0)\}_{i=0}^{n}$. Furthermore,
\begin{align}
&\bE\sbra{\prod_{x\in N_s^{(n+1)}(\cdot\mid 0)}\bE\sbra{\exp\cbra{\int\log z(x+y)\sum_{i=0}^{n}N_c^{(i)}(\d y\mid t_0)}}}
\label{}\\
=&\bE\sbra{\prod_{x\in N_s^{(n+1)}(\cdot\mid 0)}G_c^{(n)}\sbra{z(x+\cdot)\mid t_0}}
\label{}\\
=&\bE\sbra{\prod_{x\in N_s^{(n+1)}(\cdot\mid 0)}G_c^{(n)}\sbra{z_x(\cdot)\mid t_0}},
\end{align}
where in the last expression we introduced $z_x(\cdot):=z(x+\cdot)$. We recognize the p.g.fl. of the process $N_s(\cdot\mid 0)$ and from equation $\eqref{PoissonComponent}$, we have
\begin{align}
G_c^{(n+1)}[z\mid t_0]=&z(t_0)G_s\sbra{G_c^{(n)}\sbra{z_{\cdot}\mid t_0}\mid 0}
\label{}\\
=&z(t_0)\exp\cbra{\int_0^\infty \rbra{G_c^{(n)}[z_x\mid t_0]-1}h(x)\d x}.
\end{align}
Now, we want to take the limit as $n\rightarrow\infty$. First notice that for any measurable function $z$ such that $0\le z(y)\le 1$ for all $y\ge0$ and $1-z$ vanishes outside a bounded set, it holds that
\begin{align}
\lim_{n\rightarrow\infty}\prod_{i=0}^n\prod_{x\in N_c^{(i)}(\cdot\mid t_0)}z(x)=\prod_{i\ge0}\prod_{x\in N_c^{(i)}(\cdot\mid t_0)}z(x).
\end{align}
We have from the DCT,
\begin{align}
\lim_{n\rightarrow\infty}G_c^{(n)}[z\mid t_0]=\lim_{n\rightarrow\infty}\bE\sbra{\prod_{i=0}^n\prod_{x\in N_c^{(i)}(\cdot\mid t_0)}z(x)}=\bE\sbra{\prod_{i\ge0}\prod_{x\in N_c^{(i)}(\cdot\mid t_0)}z(x)}=G_c[z\mid t_0],
\end{align}
this means that $F_n$ converges to $F$ pointwise. Now, since $0\le G_c^{(n)}[z_x\mid t_0]\le 1$ for all $n\ge0$ and $\absol{G_c^{(n)}[z_x\mid t_0]-1}h(x)\le h(x)$ with $\int_0^\infty h(x)\d x<\infty$, from the DCT we get,
\begin{align}
G_c[z\mid t_0]=\lim_{n\rightarrow\infty}G_c^{(n)}[z\mid t_0]=&\lim_{n\rightarrow\infty}z(t_0)\exp\cbra{\int_0^\infty \rbra{G_c^{(n)}[z_x\mid t_0]-1}h(x)\d x}
\label{}\\
=&z(t_0)\exp\cbra{\int_0^\infty \rbra{G_c[z_x\mid t_0]-1}h(x)\d x},
\end{align} 
which concludes the proof.
}

\section{The stationary RHP}
\label{SectionStationary}
A \emph{delayed renewal process} is a renewal process in which we replace $S_0$ by a positive random variable $\hat{S}_0$ independent of $\tau_1,\tau_2,\dots$, with distribution function $F_0$ not necessarily equal to $F$. The partial sums $\hat{S}_n=\hat{S}_0+\tau_1+\dots+\tau_n$ have the associated counting process $\hat{N}_R(t)=\sum_i 1_{\cbra{\hat{S}_i\le t}}$, where the distribution of $\hat{S}_n$ for $n\ge1$ is given by $\bP(\hat{S}_n\le x)=F_0*F^{*n}(x)$ for $x\ge0$.

Under assumption \textbf{(B)}, we can obtain a stationary version of the renewal process $\hat{N}_R$ by considering a delayed renewal process with the same interarrival distributions as $N_R$, but with a suitable distribution $F_0$ for the delay $\hat{S}_0$. Such a distribution has a density (c.f. \cite[Proposition 4.2.I]{DVJ}) given as
\begin{align}
f_0=m(1-F).
\end{align}
In the stationary case, the p.g.fl. for the renewal process can be computed on the entirety of $[0,\infty)$. For this, we use the general formula \cite[(5.5.4) in Sec. V. p.146]{DVJ} to expand the p.g.fl. of the point process $\hat{N}_R$ as
\begin{align}
\label{PgflExpansion}
\hat{G}_R[z]=1+\sum_{k=1}^\infty \frac{1}{k!}\int_{\bR^{k}}(z(x_1)-1)\cdots (z(x_k)-1)M_{[k]}(\d x_1\times\cdots\times\d x_k),
\end{align}
where the \emph{factorial moment measures} $M_{[k]}(\cdot)$ for $\hat{N}_R$ are given for $A_1,\dots,A_r\in \cB(\bR)$ and nonnegative integers $k_1,\dots,k_r$ such that $k_r\ge1$ and $k_1+\dots+k_r=k$ as
\begin{align}
M_{[k]}(A_1^{(k_1)}\times\cdots\times A_r^{(k_r)})=\bE\sbra{[N(A_1)]^{[k_1]}\cdots [N(A_r)]^{[k_r]}},
\end{align}
with the factorial powers defined as
\begin{align}
n^{[k]}=\left\{\begin{array}{ll} n(n-1)\cdots(n-k+1)&(k=0,\dots, n), \\0&(k>n), \end{array}\right.
\end{align}
for any nonnegative integer $n$. We use the formula \cite[(5.4.15) in p.139]{DVJ} for the stationary renewal process and see that the factorial measures of $\hat{N}_R$ have densities on $x_1<\dots<x_k$ given by \cite[(5.4.15) in Sec. V. p.139]{DVJ}
\begin{align}
M_{[k]}(\d x_1\times\cdots\times\d x_k)= m\d x_1 \varphi(x_2-x_1)\d x_2\cdots\varphi(x_k-x_{k-1})\d x_k,
\end{align}
where we recall that $m^{-1}=\bE\sbra{\tau}$ and $\varphi=\sum_{n\ge1}f^{*n}$.
We can rewrite the integral in \eqref{PgflExpansion} using the factorial densities to obtain
\begin{align}
\hat{G}_R[z]=1+&\sum_{k=1}^\infty \frac{m}{k!}\int_0^{\infty}\int_{x_1}^\infty\cdots\int_{x_{k-1}}^{\infty}[z(x_1)-1]\cdots\nonumber
\\
&\cdots[z(x_k)-1]\d x_1 \varphi(x_2-x_1)\cdots\varphi(x_k-s_{k-1})\d x_k.
\label{}
\end{align}
Consider now a RHP in which we have replaced the centre process $N_R$ for its stationary version $\hat{N}_R$. We denote this process by
\begin{align}
\label{StationaryRHP}
\hat{N}(\cdot):=\int_0^\infty N_c(\cdot\mid t)\hat{N}_R(\d t)=\sum_{t_0\in \hat{N}_R(\cdot)}\sum_{n\ge0}N_c^{(n)}(\cdot\mid t_0),
\end{align}
whose p.g.fl. $\hat{G}[z]=\hat{G}_R[G_c[z\mid\cdot]]$ is given as
\begin{align}
\hat{G}[z]=1+&\sum_{k=1}^\infty \frac{m}{k!}\int_0^{\infty}\int_{t_1}^\infty\cdots\int_{t_{k-1}}^{\infty}(G_c[z\mid t_1]-1)\cdots\nonumber
\\
&\cdots(G_c[z\mid t_k]-1)\d t_1 \varphi(t_2-t_1)\cdots\varphi(t_k-t_{k-1})\d t_k.
\label{StationaryPgfl}
\end{align}
Since an RHP with p.g.fl. given by \eqref{StationaryPgfl} has a stationary centre process and its satellite processes satisfy $N_c(\cdot\mid t_0)\dist N_c(\cdot-t_0\mid 0)$ for $t_0\ge0$, from Vere-Jones \cite{EarthquakesVJ} we can conclude that the process is stationary and we call it the \emph{stationary renewal Hawkes process}.

\begin{Ex}
In particular, consider the renewal process $\hat{N}_R(\cdot)$ to be an homogeneous Poisson process of constant intensity $\mu>0$. In this case, the renewal density is constant $\varphi(t)=\mu$ for all $t\ge0$, and $m=\mu$. We have
\begin{align}
\hat{G}[z]=&1+\sum_{k=1}^\infty \frac{\mu^k}{k!}\int_0^\infty\cdots\int_0^\infty(G_c[z\mid s_1]-1)\cdots(G_c[z\mid s_k]-1)\d s_1\cdots\d s_k
\label{}\\
=&1+\sum_{k=1}^\infty \frac{\mu^k}{k!}\rbra{\int_0^\infty(G_c[z\mid s]-1)\d s}^k
\label{}\\
=&\sum_{k=0}^\infty \frac{1}{k!}\rbra{\int_0^\infty\mu(G_c[z\mid s]-1)\d s}^k
\label{}\\
=&\exp\cbra{\int_0^\infty\mu(G_c[z\mid s]-1)\d s}.
\end{align}
Now observe that the p.g.fl. for the satellite processes satisfies 
\begin{align}
G_c[z\mid t]=G_c[z_t\mid 0]=G_c[z(t+\cdot)\mid 0],
\end{align}
so we get
\begin{align}
\hat{G}[z]=\exp\cbra{\int_0^\infty\mu\rbra{G_c[z(s+\cdot)\mid 0]-1}\d s},
\end{align}
which is the p.g.fl. obtained by Hawkes--Oakes in \cite{HawkOak} for the linear Hawkes process. 
\end{Ex}
Finally, let us come back to the general case. We can relate the RHP with the stationary RHP in the limit through the following result.

\begin{Thm}
Let $N$ be an RHP and $\hat{N}$ be the stationary version \eqref{StationaryRHP}. Then, under assumptions \textbf{(A)} and \textbf{(B)}, 
\begin{align}
\label{ConvergenceRHP}
N(\cdot+t)\tend{d}{t\rightarrow\infty}\hat{N}(\cdot),
\end{align}
where $\cN^\sharp_{\bR}$ is equipped with the topology of vague convergence.
\end{Thm}
The proof of this Theorem follows from the convergence of the renewal process to its stationary version (see for example \cite[Sec. VI. Example 2a]{Asmussen2003}),
\begin{align}
\label{ConvergenceRenewal}
N_R(\cdot+t)\tend{d}{t\rightarrow\infty}\hat{N}_R(\cdot).
\end{align}
Heuristically, we could say that
\begin{align}
N(\cdot+t)=&\int_{\bR}N_c(\cdot+t\mid y)N_R(\d y)
\label{}\\
\dist&\int_{\bR}N_c(\cdot\mid y-t)N_R(\d y)
\label{}\\
=&\int_{\bR}N_c(\cdot\mid y)N_R(\d y+t)
\label{}\\
\tend{d}{t\rightarrow\infty}&\int_{\bR}N_c(\cdot\mid y)\hat{N}_R(\d y).
\end{align}
This result can be formalized by the convergence of the p.g.fl. Let $\Xi_c$ denote the class of continuous functions $z:\bR\rightarrow (0,1]$ s.t. $1-z$ vanishes outside a bounded set. Then, to prove \eqref{ConvergenceRHP} it is enough to show convergence of the p.g.fl. of $N$ to that of $\hat{N}$ for $z\in\Xi_c$  (c.f. \cite[Proposition 11.1.VIII]{DVJ2}), which we do as follows.
\Proof{Let $G[z]$, $\hat{G}[z]$ be the p.g.fls. of $N$, $\hat{N}$ respectively, and $G_R[z]$, $\hat{G}_R[z]$, the p.g.fls. for $N_R$ and $\hat{N}_R$ respectively. Notice that the p.g.fl. of $N(\cdot+t)$ is given by $G[z(\cdot-t)]$, then for $z\in\Xi_c$,
\begin{align}
G[z(\cdot-t)]=&G_R\sbra{G_c[z_{-t}\mid \cdot\;]}
\label{}\\
=&G_R\sbra{G_c[z\mid \cdot-t]}.
\end{align}
Note as well that $z\in\Xi_c$ implies
\begin{align}
\label{GinXi}
\tilde{z}(\cdot):= G_c[z\mid\cdot]\in \Xi_c,
\end{align}
since for $t_0\ge0$ and $z\in\Xi$,
\begin{align}
\lim_{t\rightarrow t_0}G_c[z\mid t]=&\lim_{t\rightarrow t_0}\bE\sbra{\exp\int_0^\infty\log z(s)N_c(\d s\mid t)}
\label{}\\
=&\lim_{t\rightarrow t_0}\bE\sbra{\exp\int_0^\infty\log z(s+t)N_c(\d s\mid 0)},
\end{align}
and the integrals are nonpositive, so the exponentials inside the expectation are dominated by the constant $1$. Applying the DCT once and using the continuity of the exponential function yields
\begin{align}
\lim_{t\rightarrow t_0}G_c[z\mid t]=&\bE\sbra{\exp\cbra{\lim_{t\rightarrow t_0}\int_0^\infty\log z(s+t)N_c(\d s\mid 0)}},
\label{}
\end{align}
and $N_c([0,\infty)\mid 0)$ is a.s. finite, thus applying the DCT once more and using the continuity of the logarithm and $z$, yields
\begin{align}
\lim_{t\rightarrow t_0}G_c[z\mid t]=&\bE\sbra{\exp\int_0^\infty\log z(s+t_0)N_c(\d s\mid 0)},
\label{}\\
=&\bE\sbra{\exp\int_0^\infty\log z(s)N_c(\d s\mid t_0)}
\label{}\\
=&G_c[z\mid t_0].
\end{align}
And from \eqref{ConvergenceRenewal} we have
\begin{align}
\label{ConvergencePgfl}
G_R[\tilde{z}(\cdot-t)]\tend{}{t\rightarrow\infty}\hat{G}_R[\tilde{z}]\quad\text{for } \tilde{z}\in\Xi_c.
\end{align}
In summary
\begin{align}
G[z(\cdot-t)]\tend{}{t\rightarrow\infty}\hat{G}[z]\quad\text{for all }z\in\Xi_c,
\end{align}
which proves \eqref{ConvergenceRHP}.
}
This last result implies that the p.g.fl. for the stationary RHP can be used as an approximation to the p.g.fl. of the general RHP as $t\rightarrow\infty$. It is also interesting that the limit point process is not a linear Hawkes process, since the hazard function $\mu(t)$ is not necessarily a constant unless the distribution of interarrivals $F$ is an exponential from the beginning. As an example of application of the cluster structure, we look at the first moment measure for the stationary RHP.

\begin{Prop}
Let $\hat{N}$ be a stationary RHP. Then, under \textbf{(A)} and \textbf{(B)} we have,
\begin{align}
\bE\sbra{N(B)}=\absol{B}\frac{m}{1-\alpha},\quad B\in\cB(\bR),
\end{align}
where $\absol{B}$ denotes the Lebesgue measure of $B$.
\end{Prop}
\Proof{
We compute the expectation directly using the cluster structure as
\begin{align}
\bE\sbra{N(B)}=&\bE\sbra{\sum_{t_0\in\hat{N}_R}N_c(B\mid t_0)}
\label{}\\
=&\bE\sbra{\sum_{t_0\in\hat{N}_R}N_c(B-t_0\mid 0)}.
\end{align}
We condition on the process $\hat{N}_R$ to obtain
\begin{align}
\bE\sbra{\sum_{t_0\in\hat{N}_R}N_c(B-t_0\mid 0)}=&\bE\sbra{\bE\sbra{\sum_{t_0\in\hat{N}_R}N_c(B-t_0\mid 0)\Bigg| \hat{N}_R}}
\label{}\\
=&\bE\sbra{\bE\sbra{\int_{\bR}N_c(B-t\mid 0)\hat{N}_R(\d t)\Bigg| \hat{N}_R}}
\label{}\\
=&\bE\sbra{\int_{\bR}\bE\sbra{N_c(B-t\mid 0)}\hat{N}_R(\d t)},
\label{}
\end{align}
where the conditioning on $\hat{N}_R$ disappears using the independence of the satellite processes from $\hat{N}_R$. Now we exchange the order of the integral with the expectation to write 
\begin{align}
\bE\sbra{\int_{\bR}\bE\sbra{N_c(B-t\mid 0)}\hat{N}_R(\d t)}=&\int_{\bR}\bE\sbra{N_c(B-t\mid 0)}\bE\sbra{\hat{N}_R(\d t)}.
\end{align}
From the stationarity of $\hat{N}_R$ we have that $\bE\sbra{\hat{N}_R(\d t)}=m\d t$, thus,
\begin{align}
\int_{\bR}\bE\sbra{N_c(B-t\mid 0)}\bE\sbra{\hat{N}_R(\d t)}=&m\int_{\bR}\bE\sbra{N_c(B-t\mid 0)}\d t
\label{}\\
=&m\bE\sbra{\int_{\bR}N_c(B-t\mid 0)\d t}
\label{}\\
=&m\bE\sbra{\int_{\bR}\int_{\bR}1_{B-t}(s)N_c(\d s\mid 0)\d t}
\label{}\\
=&m\bE\sbra{\int_{\bR}\int_{\bR}1_{B-s}(t)\d t N_c(\d s\mid 0)},
\end{align}
where the last equality follows from Fubini's theorem. Finally, we get,
\begin{align}
&m\bE\sbra{\int_{\bR}\int_{\bR}1_{B-s}(t)\d t N_c(\d s\mid 0)}
\label{}\\
=&m\bE\sbra{\int_{\bR}\absol{B} N_c(\d s\mid 0)}
\label{}\\
=&m\absol{B}\bE\sbra{N_c([0,\infty\mid 0)}
\label{}\\
=&m\absol{B}\frac{1}{1-\alpha}.
\end{align}
The proof is complete.
}
\section{Appendix}
\label{SectionAppendix}
We start this section by giving a proof of Lemma \ref{RenewalIntensity}.
\Proof{
From \cite[Ch.3, T7, p.61]{BremaudQueues} we know that the intensity for the process $N_R(\cdot)$ is given by
\begin{align}
\label{RHPL3E2}
\lambda_R(t)=\sum_{n\ge0}\frac{f^{(n+1)}(t-S_{n})}{1-\int_0^{t-S_{n}}f^{(n+1)}(x)\d x}1_{\cbra{S_{n}\le t<S_{{n+1}}}},
\end{align}
where for a Borel set $A$, the conditional distribution of the $n+1$-interarrival $F^{(n+1)}$ and its density $f^{(n+1)}$:
\begin{align} 
F^{(n+1)}(A):=\bP\sbra{S_{{n+1}}\in A\mid \cF_{S_{n}}}=\int_A f^{(n+1)}(x)\d x.
\end{align}
In the case of $N_R(\cdot)$, we have that $F^{(n+1)}=F$ and $f^{(n+1)}=f$ for all $n\ge0$. We can the substitute this in \eqref{RHPL3E2} to get
\begin{align}
\lambda_R(t)=&\sum_{n\ge0}\frac{f(t-S_{n})}{1-F(t-S_{n})}1_{\cbra{S_{n}\le t<S_{{n+1}}}}
\label{}\\
=&\sum_{n\ge0}\mu\rbra{t-S_{n}} 1_{\cbra{S_{n}\le t<S_{{n+1}}}}
\label{}\\
=&\sum_{n\ge0}\mu\rbra{t-S_{N_R(t)-1}} 1_{\cbra{S_{n}\le t<S_{{n+1}}}}
\label{}\\
=&\mu\rbra{t-S_{N_R(t)-1}}.
\end{align}
This completes the proof since by taking $\mu\rbra{t-S_{N_R(t-)-1}}$, we get a left-continuous modification of the intensity, and hence, the predictable version.
}

Now we present the proofs of some relevant results of Section \ref{SectionTheory}. First, we have the expression for the p.g.fl. of the independent cluster process.

\Proof[Proof of Theorem \ref{TheoremPgfl}]{
In the definition of the p.g.fl. we take expectation with respect to the law of $N(\cdot)$ given by \eqref{ClusterLaw}, and we obtain
\begin{align}
G[z]=&\bE\sbra{\exp\int\log z(x)N(\d x)}
\label{}\\
=&\int_{\cN^\sharp_{\cY}}\rbra{\bE\sbra{\exp\int_{\cX}\log z(x)N_s(\d x\mid\mu)}}\bP(N_c\in\d\mu)
\label{}\\
=&\int_{\cN^\sharp_{\cY}}\rbra{\prod_{y\in \mu(\cdot)}G_s[z\mid y]}\bP(N_c\in\d\mu)
\label{}\\
=&\bE\sbra{\prod_{y\in N_c(\cdot)}G_s[z\mid y]}
\label{}\\
=&G_c[G_s[z\mid\cdot]],
\end{align}
concluding the proof.
}

Next, we have the result from Westcott \cite{WestcottExistenceCluster} on the existence of the cluster process.

\Proof[Proof of Theorem \ref{WestcottThm}]{ Let $B\in\cB(\cX)$ bounded, then the probability generating function (p.g.f.) for the random variable $N(B)$ is obtained by evaluating the p.g.fl. of $N$ in  the function $\xi(u)=1-(1-z)1_B(u)$ where $z$ is a constant in $(0,1)$. We obtain
\begin{align}
\bE\left[z^{N(B)}\right]=\bE\left[\exp\left\{-\int_{\cY} Q_B(z\mid y)N_c(\d y)\right\}\right],
\label{}
\end{align}
where we defined $Q_B(z\mid y):=-\log\bE\left[z^{N_s(B\mid y)}\right]$. Considering now the sequence $\{z_n\}_{n\in\bN}$ given by $z_n=1-\frac{1}{n}$ for $n\ge 1$, we can compute
\begin{align}
\bP(N(B)<\infty)=&\lim_{n\rightarrow\infty}\bE\left[z_n^{N(B)}\right]
\label{}\\
=&\lim_{n\rightarrow\infty}\bE\left[\exp\left\{-\int_{\cY} Q_B(z_n\mid y)N_c(\d y)\right\}\right]
\label{}\\
=&\bE\left[\exp\left\{-\lim_{n\rightarrow\infty}\int_{\cY} Q_B(z_n\mid y)N_c(\d y)\right\}\right].
\label{}
\end{align}
The last equality is the result of applying the Monotone Convergence Theorem to the increasing sequence $\{-Q_B(z_n\mid y)\}_{n\in\bN}$ and using the continuity of the exponential function. Hence $\bP(N(B)<\infty)=1$ if and only if 
\begin{align}
\lim_{n\rightarrow\infty}\int_{\cY} Q_B(z_n\mid y)N_c(\d y)=0\quad\text{a.s.}
\label{}
\end{align}
Since for large enough $n$, $Q_B(z_n\mid y)\le Q_B(z\mid y)$ for all $y\in\cY$, and $Q_B(z_n\mid y)\downarrow 0$ as $n\rightarrow \infty$, this is equivalent, by the Dominated Convergence Theorem, to 
\begin{align}
\int_{\cY} Q_B(z\mid y) N_c(\d y)<\infty\quad\text{for some }0<z<1\quad\text{a.s.}
\label{ConditionFiniteQ}
\end{align}

Notice that if $0<a\le x\le 1$,
\begin{align}
1-x\le -\log x \le c(a) (1-x)\quad \text{for }c(a):=\frac{-\log a}{1-a}.
\label{InequalityLog}
\end{align}
Since for all $y\in\cY$, $N_s(\cX\mid 0)\ge N_s(B-y\mid 0)$, we have for all $z\in(0,1)$,
\begin{align}
\bE \sbra{z^{N_s(B\mid y)}}=\bE\sbra{z^{N_s(B-y\mid 0)}}\ge\bE \sbra{z^{N_s(\cX\mid 0)}}.
\label{}
\end{align}
We notice that because $N_s(\cX\mid 0)<\infty$ a.s.,
\begin{align}
\label{finiteclusterassump}
\bE \sbra{z^{N_s(\cX\mid 0)}}>0\quad\text{for all }z\in(0,1).
\end{align}
Thus, $\bE\sbra{ z^{N_s(B\mid y)}}$ is bounded away from zero for all $z\in(0,1)$, all $y\in\cY$ and any bounded set $B\in\cB(\cX)$, which implies that $Q_B(z\mid y)$ is finite. Moreover, \eqref{finiteclusterassump} also entails that $\tilde{c}(z):=c\left(\bE \sbra{z^{N_s(\cX\mid 0)}}\right)<\infty$, and so, we have from \eqref{InequalityLog} that
\begin{align}
1-\bE \sbra{z^{N_s(B\mid y)}}\le Q_B(z\mid y)\le \tilde{c}(z)\left(1-\bE \sbra{z^{N_s(B\mid y)}}\right),
\label{}
\end{align}
Hence, condition  \eqref{ConditionFiniteQ} holds if and only if 
\begin{align}
\int_{\cY} \left(1-\bE\sbra{ z^{N_s(B\mid y)}}\right)N_c(\d y)<\infty\quad\text{for some }0<z<1\quad\text{a.s.}
\label{ConditionFiniteP}
\end{align}
Additionally,
\begin{align}
\sum_{n=0}^\infty\bP(N_s(B\mid y)>n)z^n=&\sum_{n=0}^\infty z^n\sum_{m=n+1}^\infty\bP(N_s(B\mid y)=m)
\label{}\\
=&\sum_{m=1}^\infty \bP(N_s(B\mid y)=m)\sum_{n=0}^{m-1} z^n
\label{}\\
=&\sum_{m=1}^\infty \bP(N_s(B\mid y)=m)\frac{1-z^m}{1-z}
\label{}\\
=&\frac{1}{1-z}\left(1-\bE z^{N_s(B\mid y)}\right),
\label{}
\end{align}
since $1-z^0=0$. Hence condition  \eqref{ConditionFiniteP} holds if and only if
\begin{align}
\sum_{n=0}^\infty\left\{\int_{\cY}\bP(N_s(B\mid y)>n)N_c(\d y)\right\}z^n<\infty \quad\text{for some }0<z<1\quad\text{a.s.}
\label{ConditionFiniteE}
\end{align}
Since $\bP(N_s(B\mid y)>n)$ is decreasing in $n$, we see that condition \eqref{ConditionFiniteE} holds if and only if
\begin{align}
\int_{\cY}\bP(N_s(B\mid y)>0)N_c(\d y)<\infty\quad\text{a.s.}.
\end{align}
This concludes the proof.
}

Finally, we have the proof of the corollary to the previous theorem.

\Proof[Proof of Corollary \ref{WestcottCor}]{Notice that
\begin{align}
\bE\sbra{N_s(I\mid t)}\ge\bE\sbra{N_s(I\mid t); N_s(I\mid t)>0}\ge\bP\rbra{N_s(I\mid t)>0},
\label{RhoIneq1}
\end{align}
because $N_s(I\mid t)$ is integer-valued. From (ii), we have
\begin{align}
\bE\sbra{N_s(I \mid t)}=\bE\sbra{N_s(I-t \mid 0)}.
\label{RhoIneq2}
\end{align} 
Taking expectation in \eqref{FiniteClusterCondition} and using \eqref{RhoIneq1} and \eqref{RhoIneq2} we have
\begin{align}
\int_{\bR}\bP(N_s(I\mid t)>0)\bE\sbra{N_c(\d t)}\le&\int_{\bR}\bE\sbra{N_s(I-t \mid 0)}\bE\sbra{N_c(\d t)}
\label{}\\
=&\int_{\bR}\rbra{\int_{\bR}1_{I-t}(u)\bE\sbra{N_s(\d u \mid 0)}}\bE\sbra{N_c(\d t)}
\label{}\\
=&\int_{\bR}\rbra{\int_{\bR}1_{I-u}(t)\bE\sbra{N_c(\d t)}}\bE\sbra{N_s(\d u \mid 0)}
\label{}\\
=&\int_{\bR}\bE\sbra{N_c(I-u)}\bE\sbra{N_s(\d u \mid 0)}
\label{}\\
\le&\bE\sbra{N_s(\bR \mid 0)}\supre{}{t}\bE\sbra{N_c(I-t)}<\infty,
\end{align}
from (i) and (iii), which implies \eqref{FiniteClusterCondition}.
}

\bibliographystyle{plainreversed}

\end{document}